\documentclass[a4paper,10pt,leqno,oneside]{amsart}
\usepackage{amsmath}
\usepackage{amsfonts}
\usepackage{enumerate}
\usepackage{amsthm}
\usepackage{bm}
\usepackage{mathtools}
\usepackage{amssymb}
\usepackage[top=28truemm, bottom=25truemm, left=15truemm, right=15truemm]{geometry}
\usepackage{color}
\makeatletter \@addtoreset{equation}{section} \makeatother
\newcommand{\eref}[1]{(\ref{#1})}
\newcommand{\tref}[1]{Theorem \ref{#1}}
\newcommand{\pref}[1]{Proposition \ref{#1}}

\newcommand{\sref}[1]{Section \ref{#1}}

\theoremstyle{plain} \newtheorem{thm}{Theorem}  \newtheorem{prop}{Proposition} 
\theoremstyle{definition} \newtheorem{rem}{Remark} 
\title[Estimates of the primitive equations]{Single exponential $H^1$-upper bounds for the primitive equations}
\author[Takahito Kashiwabara]{Takahito Kashiwabara}
\address{Graduate School of Mathematical Sciences, The University of Tokyo, 3-8-1 Komaba, Meguro, 153-8914 Tokyo, Japan}
\email{tkashiwa@ms.u-tokyo.ac.jp}
\thanks{This work was supported by JSPS KAKENHI Grant-in-Aid for Scientific Research(C), Grant Number 24K06860.}

\begin{document}
\begin{abstract}
	The three dimensional primitive equations with full viscosity are considered in a horizontally periodic box $\Omega$, which are subject to either the homogeneous Neumann or Dirichlet conditions on the upper and bottom parts of the boundary.
	For a strong solution $v$ with initial data $a$, we establish \emph{a priori} bounds in $L^\infty(0, \infty; H^1(\Omega)) \cap L^2(0, \infty; \dot H^2(\Omega))$, the exponential part of which is $\exp(C \|a\|_{L^2(\Omega)}^2)$.
	This is in contrast to the upper bounds reported in the existing literature that are double exponential.
	Furthermore, the uniform-in-time estimate for the Neumann condition case, in which the Poincar\'e inequality is unavailable for $v$, seems to be new.
\end{abstract}
\maketitle

\section{Introduction}
The primitive equations are known to be governing equations for geophysical flow of large-scale atmosphere and ocean.
We refer to the review articles \cite{LiTi2018, PTZ09} for more details on physical backgrounds of the primitive equations.
In a simplified setting where variations of temperature and salinity are neglected, they are described as
\begin{equation} \label{eq: PEs}
\begin{aligned}
	\partial_t v + u \cdot \nabla v - \Delta v + \nabla_H\pi &= 0 \qquad \text{ in }\; \Omega \times (0, T), \\
	\partial_z\pi & = 0  \qquad \text{ in }\; \Omega \times (0, T), \\
	\mathrm{div}\,u & = 0  \qquad \text{ in }\; \Omega \times (0, T),
\end{aligned}
\end{equation}
where $\Omega = G \times (-h, 0)$ is a box domain in $\mathbb R^3$ with $G = (0, 1)^2$, and $h, T$ are positive constants.
Let the coordinate system be $(x_1, x_2, x_3)^\top = (x, y, z)^\top$.
The 3D velocity field is denoted by $u = (v, w)^\top = (v_1, v_2, w)^\top$, whose horizontal (in $x$- and $y$-directions) and vertical (in $z$-direction) components are $v : \Omega \to \mathbb R^2$ and $w : \Omega \to \mathbb R$, respectively.
The pressure field, which is independent of $z$, is represented by $\pi : G \to \mathbb R$, and it is enforced that $\int_G \pi \, dxdy = 0$ to ensure its uniqueness.
The 3D and 2D gradient operators are written as $\nabla = (\nabla_H, \partial_z)^\top = (\partial_x, \partial_y, \partial_z)^\top$; accordingly, the 3D and 2D  divergence operators, as well as the laplacians, are denoted by $\operatorname{div} u = \operatorname{div}_H v + \partial_z w$, $\operatorname{div}_H v = \partial_x v_1 + \partial_y v_2$, $\Delta = \operatorname{div} \nabla$, and $\Delta_H = \operatorname{div}_H \nabla_H$.
In particular, the convection term means $u \cdot \nabla v = (u \cdot \nabla) v = v \cdot \nabla_Hv + w \partial_z v$.

This is supplemented with the initial condition
\begin{equation} \label{eq: initial condition}
	v|_{t = 0} = a \qquad \text{ in }\; \Omega \times \{0\},
\end{equation}
and with one of the following sets of boundary conditions:
\begin{equation} \label{eq: Neumann BC}
\begin{alignedat}{2}
	&\partial_z v =   0, \quad w=0  &&\text{ on }\; (\Gamma_u \cup \Gamma_b) \times (0,T), \\
	&\text{$u$, $\pi$ are periodic} \qquad   && \text{ on }\; \Gamma_l \times (0,T),
\end{alignedat}
\end{equation}
which we refer to as \emph{the Neumann BC case}, and
\begin{equation} \label{eq: Dirichlet BC}
\begin{alignedat}{2}
	&\partial_z v =0, \quad w=0  &&\text{ on }\; \Gamma_u \times (0,T), \\
	&v = 0, \quad w=0 && \text{ on }\; \Gamma_b \times (0,T), \\
	&\text{$u$, $\pi$ are periodic} \qquad && \text{ on }\; \Gamma_l \times (0,T),
\end{alignedat}
\end{equation}
which we refer to as \emph{the Dirichlet BC case}.
Here, $\Gamma_u = G \times \{0\}$, $\Gamma_b = G \times \{-h\}$, and $\Gamma_l = \partial G \times [-h, 0]$ denote the upper, bottom, and lateral parts of the boundary of $\Omega$, respectively.
We could call \eref{eq: Dirichlet BC} mixed Dirichlet--Neumann boundary conditions as in \cite{GGHHK2017}; nevertheless, we employ the above terminology since a suitable reflection principle (cf.\ \cite[proof of Lemma 4.3]{GGHHK2017}) can reduce \eref{eq: Dirichlet BC} to a pure Dirichlet problem.

It is well known that the original system \eref{eq: PEs} can be reformulated, in terms of the prognostic variable $v$, as
\begin{equation} \label{eq: reformulated PEs}
\begin{alignedat}{2}
	&\partial_t v + v \cdot \nabla_H v + w \partial_z v - \Delta v + \nabla_H \pi = 0 \qquad &&\text{ in }\; \Omega \times (0,T), \\
	&w = -\int_{-h}^z \operatorname{div}_H v \, d\zeta \quad &&\text{ in }\; \Omega \times (0,T), \\
	&\operatorname{div}_H \bar v = 0 \quad &&\text{ in }\; G \times (0,T),
\end{alignedat}
\end{equation}
where $\bar v = \bar v(x, y) := \frac1h \int_{-h}^0 v(x, y, z) \, dz$ denotes the vertical average of $v$ in the last equality (it results from the 3D divergence-free condition $\operatorname{div} u = 0$ combined with the impermeability condition $w = 0$ on $\Gamma_u \cup \Gamma_b$).
Based on this representation, \cite{GMR01} proved local-in-time existence and uniqueness of a strong solution in the maximal $L^2$-regularity class, that is, $v \in L^\infty(0, T; H^1(\Omega)) \cap L^2(0, T; H^2(\Omega))$ by the Galerkin method (one can also show $v \in H^1(0, T; L^2(\Omega))$).
Later, in the cerebrated work by Cao and Titi \cite{CaTi07}, this solution was extended globally in time by establishing an \emph{a priori} estimate in the same class.
More details on earlier and recent developments in mathematical analysis of the primitive equations can be found in \cite{LiTi2018, PTZ09}.

The present work focuses on the structure of the key \emph{a priori} estimate of $v \in L^\infty(0, T; H^1(\Omega)) \cap L^2(0, T; H^2(\Omega))$, which we call \emph{the $H^1$-upper bound for the primitive equations}.
To be more specific, the $H^1$-upper bound, designated as $K_V(t)$ in \cite{CaTi07}, is triple exponential with respect to $a$ and double exponentially growing in $t$.
In fact, from equations (78), (76), (67), and (58) of \cite{CaTi07} we see that $K_V(t) \lesssim e^{K_1(t) K_z(t)}$, $K_1(t) \lesssim \|a\|_{L^2(\Omega)}^2$, $K_z(t) \lesssim e^{K_6(t)^{2/3} t}$, and $K_6(t) \lesssim e^{K_1(t)^2}$ (we consider only the principal part and assume for simplicity that the temperature equals $0$), so that
\begin{equation*}
	K_V(t) \lesssim \exp\big[ \|a\|_{L^2(\Omega)}^2 \exp \big(\exp( 2/3 \|a\|_{L^2(\Omega)}^4) t \big) \big].
\end{equation*}
In \cite{Ju07}, thanks to the uniform Gronwall inequality, a uniform-in-time $H^1$-upper bound is obtained, thus eliminating the dependence on $t$ from above.
We remark that the Poincar\'e inequality $\|v\|_{L^2(\Omega)} \le \|\nabla v\|_{L^2(\Omega)}$ is available in the setting of \cite{CaTi07, Ju07} due to the slip boundary condition $v \cdot n_{\Gamma_l} = 0$ and $\frac{\partial v}{\partial n_{\Gamma_l}} \times n_{\Gamma_l} = 0$ on $\Gamma_l$ ($n_{\Gamma_l}$ is the outer unit normal vector to $\Gamma_l$), which enables the use of the uniform Gronwall inequality.
However, the $H^1$-bound of \cite{Ju07} is still at least double exponential since its \emph{a priori} estimates follow those of \cite{CaTi07}.

For the Dirichlet BC case, \cite{KuZi07} established global well-posedness by combining the $L^6$-estimate for $v$ and $L^2$-estimate for $v_z := \partial_z v$.
One of the key arguments in \cite{KuZi07} is introduction of $\delta$ (see \cite[(2.12)]{KuZi07}) such that
\begin{equation*}
	\int_t^{t + 2\delta} \|\nabla v(s)\|_{L^2(\Omega)}^2 \, ds \le (\text{sufficiently small constant}) \qquad \forall t \in [0, T_{\max}].
\end{equation*}
However, because the way such $\delta$ depends on $T_{\max}$ or on initial data $a$ is \emph{a priori} unknown, the final $H^1$-upper bound obtained in \cite[p.\ 2748]{KuZi07}, the principal part of which reads $O(e^{\delta^{-2}})$, can grow larger than single exponential functions of $T_{\max}$ and $a$.
A different strategy was shown by \cite{HiKa2016}, which combines the $L^4$-estimate for $v$ (in more precise, for $\tilde v := v - \bar v$) and $L^2$-estimate for $v_z$.
As in \cite[p.\ 1108]{HiKa2016}, the obtained $H^1$-upper bound is double exponential and exponentially depends on $T$, whereas exponential decay is proved for the long time behavior as $t \to \infty$ by utilizing $\int_0^\infty \|v(s)\|_{H^1(\Omega)}^2 \, ds < \infty$ and a small-data-global-existence theory of the primitive equations (the Poincar\'e inequality is also essential in this result).

More refined \emph{a priori} estimates for the primitive equations including temperature are presented in \cite{CLT2016}, where only the horizontal viscosity (and diffusivity) is considered.
Namely, a more difficult problem than \eref{eq: PEs} is treated there (see also \cite{CLT2017, HSW2021}).
However, as far as the structure of $H^1$-upper bounds is concerned, the bounds obtained in \cite{CLT2016, CLT2017, HSW2021} are at least double exponential because of the use of the logarithmic Sobolev and Gronwall inequalities (see \cite[Lemmas 2.2 and 2.3]{CLT2016}).
Recently, long time behavior of solutions with some smallness assumptions on initial data is examined in \cite{TanLiu2026}.

In this paper, we consider the primitive equations \eref{eq: PEs} with the full viscosity and establish $H^1$-upper bounds that are only single exponential with respect to $a$ and are independent of $t$ for both of the Neumann and Dirichlet BC cases.
Namely, the following is our main result.

\begin{thm} \label{main thm}
	Suppose that the initial data satisfies $a \in H^1(\Omega)$ and $a_z: = \partial_z a \in L^3(\Omega)$, as well as the compatibility conditions: $\operatorname{div}_H \bar a = 0$, the horizontal periodicity, and $a = 0$ on $\Gamma_b$ in case of \eref{eq: Dirichlet BC}.
	Then the solution of \eref{eq: PEs}--\eref{eq: initial condition} with \eref{eq: Neumann BC} or \eref{eq: Dirichlet BC} satisfies 
	\begin{equation*}
		\|v(t)\|_{H^1(\Omega)}^2 + \int_0^t \|\nabla^2 v(s)\|_{L^2(\Omega)}^2 \, ds \le B(a) e^{C \|a\|_{L^2(\Omega)}^2} \quad \forall t \in [0, T],
	\end{equation*}
	where $B(a)$ is a polynomial function of $\|a\|_{H^1(\Omega)}$ and $\|a_z\|_{L^3(\Omega)}$, and $C$ is a positive constant dependent only on $\Omega$.
\end{thm}
\begin{rem}
	(i) Since the above $H^1$-bound is independent of $T$, one can even have
	\begin{equation*}
		\|v(t)\|_{H^1(\Omega)}^2 + \int_0^\infty \|\nabla^2 v(s)\|_{L^2(\Omega)}^2 \, ds \le B(a) e^{C \|a\|_{L^2(\Omega)}^2} \quad \forall t \in [0, \infty).
	\end{equation*}
	
	(ii) $B(a)$ actually depends on $\|a\|_{L^4(\Omega)}$ and $\|a\|_{L^6(\Omega)}$ (see e.g.\ \eref{eq: L4 estimate for v in Drichlet BC case} below), which can be bounded by $\|a\|_{H^1(\Omega)}$ in virtue of the Sobolev embedding theorem.
	
	(iii) For the Neumann BC case, since $\int_\Omega a \, dxdydz \neq 0$ in general, we do not expect for $v$ to satisfy the Poincar\'{e} inequality as in \cite[(3.174) and (3.279)]{PTZ09}.
	To the best of our knowledge, a uniform-in-time $H^1$-upper bound for such a situation was not shown in the existing literature.
	
	(iv) For the Dirichlet BC case, if the Poincar\'e inequality is taken into account, an exponential decaying bound as $t \to \infty$ could be derived.
	Nevertheless it is not directly addressed in this paper.
\end{rem}

As a building block to prove the theorem above, we first obtain in \sref{sec: proof of main prop} the following result for the Neumann BC case.
\begin{prop} \label{main prop}
	If $a \in H^1(\Omega) \cap L^8(\Omega)$, $a_z \in L^4(\Omega)$, $\operatorname{div}_H \bar a = 0$, and the horizontal periodicity of $a$ hold, then the solution of \eref{eq: PEs}--\eref{eq: initial condition} with \eref{eq: Neumann BC} satisfies 
	\begin{equation*}
		\|v(t)\|_{H^1(\Omega)} \le B(a, t) e^{C (\sqrt{t} + 1) \|a\|_{L^2(\Omega)}^2} \quad \forall t \in [0, T],
	\end{equation*}
	where $B(a, t)$ is a polynomial functions of $\|a\|_{H^1(\Omega)}, \|a\|_{L^8(\Omega)}, \|a_z\|_{L^4(\Omega)}, \sqrt{t}$, and $C$ is a positive constant dependent only on $\Omega$.
\end{prop}

The proof of \pref{main prop}, in which we utilize the $L^4(\Omega)$- and $L^8(\Omega)$-estimates for $v$ as well as $L^2(\Omega)$- and $L^4(\Omega)$-estimates for $v_z$, is more straightforward than that of \tref{main thm}.
In fact, the dissipative parts associated with those estimates are nicely consistent in dealing with the ``most difficult nonlinearity'' of the primitive equations, i.e., $w v_z$ (see \eref{eq1: treatment of wv_z} and \eref{eq2: treatment of wv_z} below).
However, there are two drawbacks with this approach: 1) the technique is not directly applicable to the Dirichlet BC case, and 2) the bound depends on $t$.

In \sref{sec1: proof of main thm} we overcome the first difficulty and prove \tref{main thm} for the Dirichlet BC case.
The idea is to make use of the $L^6(\Omega)$-estimate for $v$ and $L^3(\Omega)$-estimate for $v_z$ (instead of $L^8$- and $L^4$-estimates).
Employing a similar strategy, in \sref{sec2: proof of main thm} we eliminate time dependency from the $H^1$-upper bound and prove \tref{main thm} for the Neumann BC case.
Since the main issue here is treatment of ``lower order terms'' (see e.g.\ \eref{eq: treatment of lower-order term}), a Poincar\'e-type inequality is a key ingredient for this observation.
Thus we are led to consider $L^4(\Omega)$- and $L^6(\Omega)$-estimates for $\tilde v = v - \bar v$ in place of $v$, which satisfies $\|\tilde v\|_{L^2(\Omega)} \le h \|v_z\|_{L^2(\Omega)}$ (cf.\ \cite[(6.4)]{GGHHK2020}).

We will exploit the standard notation of the Lebesgue spaces $L^p(D)$ and the Sobolev spaces $W^{s, p}(D)$ in a Lipschitz domain $D$.
Norms of a function $f$ in isotropic spaces are denoted as
\begin{equation*}
	\|f\|_{L^p_z L^q_{x, y}} = \Big( \int_{-h}^0 \|f(z)\|_{L^q(G)}^p \, dz \Big)^{1/p} = \Big( \int_{-h}^0 \Big( \int_G |f(x, y, z)|^q \, dxdy \Big)^{p/q} \, dz \Big)^{1/p},
\end{equation*}
with appropriate modifications for the cases $p, q = \infty$.
For integration in $\Omega$ and $G$, we will omit the symbols of volume elements, i.e., $dx dy dz$ and $dx dy$.
The $L^2$-inner products in $\Omega$ and $G$ are denoted by $(\cdot, \cdot)_\Omega$ and $(\cdot, \cdot)_G$, respectively.

Finally, we recall the Gronwall inequality of the following form: if nonnegative functions $y(t), z(t), \alpha(t)$, and $\beta(t)$ satisfy $\frac{dy}{dt}(t) + z(t) \le \alpha(t) y(t) + \beta(t)$ for $t > 0$, then for all $t \ge 0$ we have
\begin{align*}
	y(t) &\le \Big( y(0) + \int_0^t \beta(s) \, ds \Big) e^{\int_0^t \alpha(s) \, ds}, \\ 
	y(t) + \int_0^t z(s) \, ds &\le \Big( y(0) + \int_0^t \beta(s) \, ds \Big) \Big(1 + \int_0^t \alpha(s) \, ds \, e^{\int_0^t \alpha(s) \, ds} \Big).
\end{align*}

\section{Proof of \pref{main prop}} \label{sec: proof of main prop}
To prove \pref{main prop} (and \tref{main thm} as well) we may assume that a sufficiently smooth solution for \eref{eq: PEs}--\eref{eq: initial condition} with \eref{eq: Neumann BC} or \eref{eq: Dirichlet BC} uniquely exists; in fact, this is already known by \cite[Theorems 3.1 and 3.3]{GGHHK2020}.
In the following, $C$ denotes a generic positive constant which may depend only on $\Omega$ unless otherwise stated.

\subsection{Step 1: estimates for $v \in L^\infty_t L^q(\Omega)$} \label{subsec2.1}
Although this step is essentially known (see \cite[Proposition 3.1(iii)]{CLT2016}), we present the calculation for completeness.
First we multiply \eref{eq: reformulated PEs} by $v$ and integrate over $\Omega \times [0, t]$ to obtain
\begin{equation} \label{eq: energy equality}
	\frac12 \|v(t)\|_{L^2(\Omega)}^2 + \int_0^t \|\nabla v(s)\|_{L^2(\Omega)}^2 \, ds = \frac12 \|a\|_{L^2(\Omega)}^2.
\end{equation}

Next, to get an expression of $\pi$ in terms of $v$, we apply $\frac{1}{h} \operatorname{div}_H \int_{-h}^0 (\cdot) \, dz$ to \eref{eq: reformulated PEs}$_1$.
It follows that, since $\int_{-h}^0 w v_z \, dz = \int_{-h}^0 (\operatorname{div}_H v) v \, dz$,
\begin{equation} \label{eq: pressure Poisson for Neumann BC case}
	-\Delta_H \pi = \frac1h \operatorname{div}_H \Big( \int_{-h}^0 \operatorname{div}_H (v \otimes v) \, dz \Big), \quad
	\nabla_H \pi = \frac1h \big[ \nabla_H (-\Delta_H)^{-1} \operatorname{div}_H \big] \int_{-h}^0 \operatorname{div}_H (v \otimes v) \, dz,
\end{equation}
where $\operatorname{div}_H (v \otimes v)$ means the vector whose $i$-th component is $\sum_{j = 1}^2 \partial_{x_j} (v_i v_j)$ for $i = 1, 2$.
We then apply the Calderon--Zygmund inequality (cf.\ \cite[Lemma 4.1]{GGHHK2017}) and the Minkowski inequality to deduce that
\begin{equation} \label{eq: pressure estimate in Neumann BC case}
	\|\nabla_H \pi\|_{L^q(G)} \le C \Big\| \int_{-h}^0 \operatorname{div}_H (v \otimes v) \, dz \Big\|_{L^q(G)} \le C \int_{-h}^0 \big\| |v| |\nabla_H v| \big\|_{L^q(G)} \, dz \quad (q \in (1, \infty)),
\end{equation}
where the constant $C$ may depend on $q$.

Let us multiply \eref{eq: reformulated PEs}$_1$ by $|v|^{q-2} v \, (q > 2)$ and use the cancelation property of the convection term, that is, $(u \cdot \nabla v, |v|^{q-2} v)_\Omega = 0$ for the divergence-free 3D vector $u$ with the vanishing (or periodic) normal components on the boundary.
It follows from integration by parts that
\begin{equation} \label{eq: start of Lq estimate}
	\frac1q \frac{d}{dt} \|v\|_{L^q(\Omega)}^q + \frac{4(q-2)}{q^2} \big\| \nabla |v|^{\frac{q}2} \big\|_{L^2(\Omega)}^2 + \big\| |v|^{\frac{q}2 - 1} |\nabla v| \big\|_{L^2(\Omega)}^2
		= - \Big( \nabla_H \pi, \int_{-h}^0 |v|^{q-2} v \, dz \Big)_G.
\end{equation}
Substituting \eref{eq: pressure estimate in Neumann BC case} and using the H\"older inequality with $\frac1q + \frac12 + \frac{q-2}{2q} = 1$, we see that the right-hand side above is bounded by
\begin{equation} \label{eq1: proof of main prop}
	\|\nabla_H \pi\|_{L^{ \frac{2q}{q+2} }(G)} \Big\| \int_{-h}^0 |v|^{q-2} v \, dz \Big\|_{L^{\frac{2q}{q-2}}(G)}
		\le \int_{-h}^0 \|v\|_{L^q(G)} \|\nabla_H v\|_{L^2(G)} \, dz \int_{-h}^0 \big\| |v|^{q-1} \big\|_{L^{\frac{2q}{q-2}}(G)} \, dz.
\end{equation}
The first integral is bounded by $\|1\|_{L^{\frac{2q}{q-2}}_z} \|v\|_{L^q(\Omega)} \|\nabla_H v\|_{L^2(\Omega)}$.
To estimate the second one, by the embedding $W^{1, \frac{q}{q-1}}(G) \hookrightarrow L^{ \frac{2q}{q-2} }(G)$ as well as the H\"older inequality with $\frac{q-2}{2q} + \frac12 = \frac{q-1}{q}$, we have
\begin{align*}
	\big\| |v|^{q-1} \big\|_{L^{\frac{2q}{q-2}}(G)} &\le C ( \big\| |v|^{q-1} \big\|_{L^{\frac{q}{q-1}}(G)} + \big\| \nabla_H |v|^{q-1} \big\|_{L^{\frac{q}{q-1}}(G)}) \\
		&\le C ( \|v\|_{L^q(G)}^{q-1} + \big\| |v|^{\frac{q}2 - 1} \nabla_H |v|^{\frac{q}2} \big\|_{L^{\frac{q}{q-1}}(G)} ) \\
		&\le C ( \|v\|_{L^q(G)}^{q-1} + \|v\|_{L^q(G)}^{\frac{q}2 - 1} \big\| \nabla_H |v|^{\frac{q}2} \big\|_{L^2(G)} ).
\end{align*}
Consequently,
\begin{equation*}
	\int_{-h}^0 \big\| |v|^{q-1} \big\|_{L^{\frac{2q}{q-2}}(G)} \, dz \le C \|1\|_{L^{q}_z} (\|v\|_{L^q(\Omega)}^{q-1} + \|v\|_{L^q(\Omega)}^{\frac{q}2 - 1} \big\| \nabla_H |v|^{\frac{q}2} \big\|_{L^2(\Omega)}).
\end{equation*}
Since the right-hand side of \eref{eq1: proof of main prop} is bounded by
\begin{equation*}
	C (\|\nabla_H v\|_{L^2(\Omega)} + \|\nabla_H v\|_{L^2(\Omega)}^2) \|v\|_{L^q(\Omega)}^q + \frac{q - 2}{q^2} \big\| \nabla_H |v|^{\frac{q}2} \big\|_{L^2(\Omega)}^2,
\end{equation*}
we obtain, for all $t > 0$,
\begin{equation} \label{eq: lower-order term in Lq estimate}
	\frac1q \frac{d}{dt} \|v\|_{L^q(\Omega)}^q + \frac{3(q-2)}{q^2} \big\| \nabla |v|^{\frac{q}2} \big\|_{L^2(\Omega)}^2 + \big\| |v|^{\frac{q}2 - 1} |\nabla v| \big\|_{L^2(\Omega)}^2 \le C (\|\nabla_H v\|_{L^2(\Omega)}^2 + \|\nabla_H v\|_{L^2(\Omega)}) \|v\|_{L^q(\Omega)}^q.
\end{equation}
Now the Gronwall inequality and $\int_0^t \|\nabla v(s)\|_{L^2(\Omega)}^2 \, ds \le \frac{\|a\|_{L^2(\Omega)}^2}{2}$ (recall \eref{eq: energy equality}) conclude, for all $t \ge 0$ and $q \in (2, \infty)$,
\begin{equation} \label{eq: Lq bound of v for Neumann BC case}
\begin{aligned}
	&\|v(t)\|_{L^q(\Omega)}^q + \int_0^t \big\| |v|^{\frac{q}2 - 1} |\nabla v| \big\|_{L^2(\Omega)}^2 \, ds \\
	\le \;& C \big( \|a\|_{L^q(\Omega)}^q + (\sqrt t + 1)\|a\|_{L^2(\Omega)}^2 \big)
		\big(1 + (\sqrt t + 1) \|a\|_{L^2(\Omega)}^2 e^{C (\sqrt{t} + 1) \|a\|_{L^2(\Omega)}^2} \big)=: A_1(q, a, t),
\end{aligned}
\end{equation}
provided that the initial data satisfy $a \in L^q(\Omega)$.
In particular, we will make use of these estimates for $q = 4$ and $q = 8$ below.

\subsection{Step 2: estimates for $v_z := \partial_z v \in L^\infty_t L^r(\Omega) \, (r = 2, 4)$}
Although we will restrict to the cases $r = 2$ and $r = 4$ in the end, we first multiply \eref{eq: reformulated PEs}$_1$ by $-\partial_z( |v_z|^{r-2} v_z )$ and perform integration by parts in $\Omega$ with general $r \ge 2$.
Observing that no boundary terms appear and using the cancelation property for the convection term, we obtain
\begin{equation} \label{eq2: proof of main prop}
\begin{aligned}
	\frac1r \frac{d}{dt} \|v_z\|_{L^r(\Omega)}^r + \frac{4(r-2)}{r^2} \big\| \nabla |v_z|^{\frac{r}2} \big\|_{L^2(\Omega)}^2 + \big\| |v_z|^{\frac{r}2 - 1} |\nabla v_z| \big\|_{L^2(\Omega)}^2 &= -(u_z \cdot \nabla v, |v_z|^{r-2} v_z)_\Omega \\
		&= -(v_z \cdot \nabla_H v, |v_z|^{r-2} v_z)_\Omega + (\operatorname{div}_H v, |v_z|^r)_\Omega.
\end{aligned}
\end{equation}
By further integration by parts with respect to horizontal derivatives, the right-hand side terms are bounded as
\begin{equation} \label{eq1: treatment of wv_z}
\begin{aligned}
	C \int_\Omega |v| |v_z|^{r-1} |\nabla_H v_z| &\le C \big\| |v| |v_z|^{\frac{r}2} \big\|_{L^2(\Omega)} \big\| |v_z|^{\frac{r}2 - 1}|\nabla_H v_z| \big\|_{L^2(\Omega)} \\
		&\le \frac1{10} \big\| |v|^{\frac{r}2 - 1} |\nabla v_z| \big\|_{L^2(\Omega)}^2 + C \big\| |v| |v_z|^{\frac{r}2} \big\|_{L^2(\Omega)}^2,
\end{aligned}
\end{equation}
where the first term can be absorbed into the left-hand side.

To treat the second term, we consider either $r = 2$ or $r = 4$.
If $r = 2$, we find immediately from \eref{eq: Lq bound of v for Neumann BC case} for $q = 4$ that $\int_0^t \big\| |v| |v_z| \big\|_{L^2(\Omega)}^2 \, ds \le A_1(4, a, t)$.
Therefore, just integration of \eref{eq2: proof of main prop} with respect to $t$ gives us, if $a_z \in L^2(\Omega)$,
\begin{equation} \label{eq4: proof of main prop}
	\frac12 \|v_z(t)\|_{L^2(\Omega)}^2 + \frac9{10} \int_0^t \|\nabla v_z\|_{L^2(\Omega)}^2 \, ds
		\le \frac12 \|a_z\|_{L^2(\Omega)}^2 + C A_1(4, a, t).
\end{equation}

Next, for $r = 4$ we observe from integration by parts that
\begin{equation} \label{eq2: treatment of wv_z}
\begin{aligned}
	\big\| |v| |v_z|^2 \big\|_{L^2(\Omega)}^2 &= \int_\Omega |v|^2 (v_z \cdot v_z) |v_z|^2 \\
		&= \underbrace{ - \int_\Omega \partial_z (|v|^2) \cdot \frac12 \partial_z (|v|^2) |v_z^2| }_{ \le 0 } - \int_\Omega |v|^2 (v \cdot \partial_z v_z) |v_z|^2 - \int_\Omega |v|^2 (v \cdot v_z) 2 (v_z \cdot \partial_z v_z) \\
		&\le 3 (|v|^3 |v_z|, |v_z| |\partial_z v_z|)_\Omega \\
		&\le \frac1{10} \big\| |v_z| |\partial_z v_z| \big\|_{L^2(\Omega)}^2 + C \big\| |v|^3 |v_z| \big\|_{L^2(\Omega)}^2.
\end{aligned}
\end{equation}
Therefore, if $a_z \in L^4(\Omega)$, then \eref{eq: Lq bound of v for Neumann BC case} for $q = 8$ leads to, for all $t \ge 0$,
\begin{equation} \label{eq: L4 bound of vz for Neumann BC case}
\begin{aligned}
	\frac14 \|v_z(t)\|_{L^4(\Omega)}^4 + \int_0^t \Big( \frac12 \big\| \nabla |v_z|^2 \big\|_{L^2(\Omega)}^2 + \frac8{10} \big\| |v_z| |\nabla v_z| \big\|_{L^2(\Omega)}^2 \Big) \, ds &\le \frac14 \|a_z\|_{L^4(\Omega)}^4 + C A_1(8, a, t) =: A_2(a, t).
\end{aligned}
\end{equation}

\subsection{Step 3: estimate for $\nabla_H v \in L^\infty_t L^2(\Omega)$} \label{subsec2.3}
Multiply \eref{eq: reformulated PEs}$_1$ by $-\Delta_H v$ and integrate by parts over $\Omega$ to obtain
\begin{equation} \label{eq3: proof of main prop}
\begin{aligned}
	\frac12 \frac{d}{dt} \|\nabla_H v\|_{L^2(\Omega)}^2 + \|\nabla \nabla_H v\|_{L^2(\Omega)}^2
		&= (v \cdot \nabla_H v, \Delta_H v)_\Omega + (w v_z, \Delta_H v)_\Omega \\
		&\le \frac1{10} \|\nabla_H^2 v\|_{L^2(\Omega)}^2 + C\big\| |v| |\nabla_H v| \big\|_{L^2(\Omega)}^2 + \|w v_z\|_{L^2(\Omega)} \|\Delta_H v\|_{L^2(\Omega)},
\end{aligned}
\end{equation}
where $\nabla_H^2 v$ means the Hessian tensor of rank $3$ for a vector-valued function $v$.
We apply an anisotropic estimate (cf.\ \cite[p.\ 1108]{HiKa2016}) to handle the last term as
\begin{equation} \label{eq: estimate of wvz with L4}
\begin{aligned}
	\|w v_z\|_{L^2(\Omega)} &\le C \|w\|_{L^\infty_z L^4_{xy}} \|v_z\|_{L^4(\Omega)} \le C \|\operatorname{div}_H v\|_{L^2_z L^4_{xy}} \|v_z\|_{L^4(\Omega)} \\
		&\le C \|\nabla_H v\|_{L^2(\Omega)}^{1/2} \|\nabla_H v\|_{H^1(\Omega)}^{1/2} \|v_z\|_{L^4(\Omega)},
\end{aligned}
\end{equation}
so that
\begin{align*}
	\|w v_z\|_{L^2(\Omega)} \|\Delta_H v\|_{L^2(\Omega)} &\le C \|\nabla_H v\|_{L^2(\Omega)}^{1/2} \|v_z\|_{L^4(\Omega)} (\|\nabla \nabla_H v\|_{L^2(\Omega)}^{3/2} + \|\nabla_H v\|_{L^2(\Omega)}^{1/2} \|\Delta_H v\|_{L^2(\Omega)}) \\
		&\le \frac1{10} \|\nabla \nabla_H v\|_{L^2(\Omega)}^2 + C \|\nabla_H v\|_{L^2(\Omega)}^2 (\|v_z\|_{L^4(\Omega)}^4 + 1).
\end{align*}
Therefore, just integration of \eref{eq3: proof of main prop} with respect to $t$, combined with \eref{eq: Lq bound of v for Neumann BC case} for $q = 4$ and with \eref{eq: L4 bound of vz for Neumann BC case} yields, for all $t \ge 0$,
\begin{align*}
	\frac12 \|\nabla_H v(t)\|_{L^2(\Omega)}^2 + \frac8{10} \int_0^t \|\nabla \nabla_H v\|_{L^2(\Omega)}^2 \, ds
		&\le \frac12 \|\nabla_H a\|_{L^2(\Omega)}^2
			+ C A_1(4, a, t) + C \|a\|_{L^2(\Omega)}^2 (A_2(a, t) + 1).
\end{align*}
This together with \eref{eq4: proof of main prop} concludes the desired upper bound for $v \in L^\infty(0, T; H^1(\Omega)) \cap L^2(0, T; H^2(\Omega))$.

\section{Proof of \tref{main thm} for Dirichlet BC case} \label{sec1: proof of main thm}
Let us turn our attention to the boundary condition \eref{eq: Dirichlet BC}, in which case the pressure $\pi$ involves an additional component.
In fact, the Poisson equation in \eref{eq: pressure Poisson for Neumann BC case} is replaced as
\begin{equation*}
	-\Delta_H \pi = \frac1h \operatorname{div}_H \Big( \int_{-h}^0 \operatorname{div}_H (v \otimes v) \, dz \Big) - \frac1h \operatorname{div}_H (v_z|_{\Gamma_b}),
\end{equation*}
that is,
\begin{equation*}
	\nabla_H \pi = \frac1h \big[ \nabla_H (-\Delta_H)^{-1} \operatorname{div}_H \big] \Big( \int_{-h}^0 \operatorname{div}_H (v \otimes v) \, dz \Big)
			- \frac1h \big[ \nabla_H (-\Delta_H)^{-1} \operatorname{div}_H \big] (v_z|_{\Gamma_b})
		=: \nabla_H \pi_1 + \nabla_H \pi_2.
\end{equation*}
The Calderon--Zygmund inequality implies that, for $q \in (0, \infty)$,
\begin{equation*}
	\|\nabla_H \pi_1\|_{L^q(G)} \le C \int_{-h}^0 \big\| |v| |\nabla_H v| \big\|_{L^q(G)} \, dz, \qquad
	\|\nabla_H \pi_2\|_{L^q(G)} \le C \|v_z\|_{L^q(\Gamma_b)}.
\end{equation*}

On the one hand, since the appearance of $\pi_2$---especially $(\nabla_H \pi_2, |v|^{q-2} v)_\Omega$---prevents the $L^q(\Omega)$-estimate for $v$ from being closed solely by itself, we need to combine it with $L^r(\Omega)$-estimate for $v_z$.
On the other hand, in the latter estimate we encounter $(\nabla_H \pi_1, |v_z|^{r-2} v_z)_{\Gamma_b}$, which need be taken care by the former one.
Below we show that choosing the two combinations $(q, r) = (4, 2), (6, 3)$ enables us to resolve this delicate interdependence between the two estimates, which was absent in the Neumann BC case.

\subsection{Step 1: estimates for $v \in L^\infty_t L^4(\Omega)$ and $v_z \in L^\infty_t L^2(\Omega)$}
The energy equality \eref{eq: energy equality} remains the same, and we also have $\int_0^t \|v\|_{H^1(\Omega)}^2 \, ds \le C \|a\|_{L^2(\Omega)}^2 \, (\forall t \ge 0)$ as a result of the Poincar\'e inequality.
Addressing $\pi_1$ in the same way as before, in place of \eref{eq: start of Lq estimate} for $q = 4$ we obtain
\begin{equation} \label{eq2: proof of main thm for Dirichlet BC case}
	\frac14 \frac{d}{dt} \|v\|_{L^4(\Omega)}^4 + \frac38 \big\| \nabla |v|^2 \big\|_{L^2(\Omega)}^2 + \big\| |v| |\nabla v| \big\|_{L^2(\Omega)}^2 \le C (\|\nabla_H v\|_{L^2(\Omega)}^2 + \|\nabla_H v\|_{L^2(\Omega)}) \|v\|_{L^4(\Omega)}^4 - \Big( \nabla_H \pi_2, \int_{-h}^0 |v|^2 v \, dz \Big)_G.
\end{equation}
By the estimate of $\nabla_H \pi_2$ mentioned above and a trace inequality, we estimate the last term as
\begin{equation} \label{eq1: proof of main thm for Dirichlet BC case}
	C \|v_z\|_{L^2(\Gamma_b)} \int_{-h}^0 \big\| |v|^3 \big\|_{L^2(G)} dz \le C \|v_z\|_{L^2(\Omega)}^{1/2} \|\nabla v_z\|_{L^2(\Omega)}^{1/2}
		\int_{-h}^0 \big\| |v|^2 \big\|_{L^3(G)}^{3/2} \, dz.
\end{equation}
Since $H^{1/3}(G) \hookrightarrow L^3(G)$, we have
\begin{align*}
	\big\| |v|^2 \big\|_{L^3(G)}^{3/2} &\le C \big\| |v|^2 \big\|_{L^2(G)} \big( \big\| |v|^2 \big\|_{L^2(G)} + \big\| \nabla_H |v|^2 \big\|_{L^2(G)} \big)^{1/2} \\
		&\le C \|v\|_{L^4(G)}^2 \big( \|v\|_{L^4(G)} + \big\| \nabla_H |v|^2 \big\|_{L^2(G)}^{1/2} \big).
\end{align*}
Hence the right-hand side of \eref{eq1: proof of main thm for Dirichlet BC case} is bounded by
\begin{align*}
	&C \|v_z\|_{L^2(\Omega)}^{1/2} \|\nabla v_z\|_{L^2(\Omega)}^{1/2} \cdot \|1\|_{L^4_z}
		C \|v\|_{L^4(\Omega)}^2 \big( \|v\|_{L^4(\Omega)} + \big\| \nabla_H |v|^2 \big\|_{L^2(\Omega)}^{1/2} \big) \\
	\le \; & \frac1{10} \|\nabla v_z\|_{L^2(\Omega)}^2 + \frac1{8} \big\| \nabla_H |v|^2 \big\|_{L^2(\Omega)}^2
		+ C\|v_z\|_{L^2(\Omega)} \|v\|_{L^4(\Omega)}^4 + C\|v_z\|_{L^2(\Omega)}^{2/3} \|v\|_{L^4(\Omega)}^4.
\end{align*}
This together with \eref{eq2: proof of main thm for Dirichlet BC case}, with $\|v\|_{L^4(\Omega)} \le C \|v\|_{H^1(\Omega)}$, and with $\alpha^2 \beta^3 \le \frac34 \alpha^2 \beta^4 + \frac14 \alpha^2$ leads to
\begin{equation} \label{end of L4 estimate for Dirichlet BC case}
	\frac14 \frac{d}{dt} \|v\|_{L^4(\Omega)}^4 + \frac28 \big\| \nabla |v|^2 \big\|_{L^2(\Omega)}^2 + \big\| |v| |\nabla v| \big\|_{L^2(\Omega)}^2 \le C \|v\|_{H^1(\Omega)}^2 (\|v\|_{L^4(\Omega)}^4 + 1) + \frac1{10} \|\nabla v_z\|_{L^2(\Omega)}^2.
\end{equation}

Next, in place of \eref{eq2: proof of main prop} for $r = 2$ we have
\begin{equation} \label{eq3: proof of main thm for Dirichlet BC case}
	\frac12 \frac{d}{dt} \|v_z\|_{L^2(\Omega)}^2 + \frac9{10} \|\nabla v_z\|_{L^2(\Omega)}^2 \le C \big\| |v| |v_z| \big\|_{L^2(\Omega)}^2 + \frac1h (\nabla_H \pi, v_z|_{\Gamma_b})_G,
\end{equation}
where we see that
\begin{gather*}
	|(\nabla_H \pi_1, v_z|_{\Gamma_b})_G| \le C \big\| |v| |\nabla_H v| \big\|_{L^2(\Omega)} \|v_z\|_{L^2(\Gamma_b)}
		\le \frac1{10} \big\| |v| |\nabla_H v| \big\|_{L^2(\Omega)}^2 + \frac1{10} \|\nabla v_z\|_{L^2(\Omega)}^2 + C \|v_z\|_{L^2(\Omega)}^2, \\
	|(\nabla_H \pi_2, v_z|_{\Gamma_b})_G| \le C \|v_z\|_{L^2(\Gamma_b)}^2 \le \frac1{10} \|\nabla v_z\|_{L^2(\Omega)}^2 + C \|v_z\|_{L^2(\Omega)}^2.
\end{gather*}
Now, addition of \eref{end of L4 estimate for Dirichlet BC case} and \eref{eq3: proof of main thm for Dirichlet BC case} deduces
\begin{align*}
	\frac{d}{dt} \Big( \frac14 \|v\|_{L^4(\Omega)}^4 + \frac12 \|v_z\|_{L^2(\Omega)}^2 \Big)
		+ \frac14 \big\| \nabla |v|^2 \big\|_{L^2(\Omega)}^2 + \frac9{10} \big\| |v| |\nabla v| \big\|_{L^2(\Omega)}^2 + \frac7{10} \|\nabla v_z\|_{L^2(\Omega)}^2
	\le C \|v\|_{H^1(\Omega)}^2 (\|v\|_{L^4(\Omega)}^4 + 1).
\end{align*}
Thus we conclude from the Gronwall inequality that, for all $t \ge 0$,
\begin{equation} \label{eq: L4 estimate for v in Drichlet BC case}
\begin{aligned}
	&\|v(t)\|_{L^4(\Omega)}^4 + \|v_z(t)\|_{L^2(\Omega)}^2 + \int_0^t \Big( \big\| |v| |\nabla v| \big\|_{L^2(\Omega)}^2 + \|\nabla v_z\|_{L^2(\Omega)}^2 \Big) \, ds \\
	\le \; & C \big( \|a\|_{L^4(\Omega)}^4 + \|a_z\|_{L^2(\Omega)}^2 + \|a\|_{L^2(\Omega)}^2 \big)
		\big( 1 + \|a\|_{L^2(\Omega)}^2 e^{C \|a\|_{L^2(\Omega)}^2} \big) =: B_1(a).
\end{aligned}
\end{equation}

\subsection{Step 2: estimates for $v \in L^\infty_t L^6(\Omega)$ and $v_z \in L^\infty_t L^3(\Omega)$}
In place of \eref{eq: start of Lq estimate} for $q = 6$ we have
\begin{equation} \label{eq4: proof of main thm for Dirichlet BC case}
	\frac16 \frac{d}{dt} \|v\|_{L^6(\Omega)}^6 + \frac49 \|\nabla |v|^3\|_{L^2(\Omega)}^2 + \big\| |v|^2 |\nabla v| \big\|_{L^2(\Omega)}^2 \le C \|v\|_{H^1(\Omega)}^2 (\|v\|_{L^6(\Omega)}^6 + 1)
		- \Big( \nabla_H \pi_2, \int_{-h}^0 |v|^4 v \, dz \Big)_G.
\end{equation}
In virtue of the trace inequality $\|v_z\|_{L^6(\Gamma_b)} = \big\| |v_z|^{3/2} \big\|_{L^4(\Gamma)}^{2/3} \le C \big\| \nabla |v_z|^{3/2} \big\|_{L^2(\Omega)}^{2/3}$, the last term is estimated, as a result of the Young and H\"older inequalities with $\frac16 + \frac56 = 1$ and the embedding $H^1(\Omega) \hookrightarrow L^6(\Omega)$, by
\begin{align*}
	\|\nabla_H \pi_2\|_{L^6(G)} \int_{-h}^0 \big\| |v|^5 \big\|_{L^{6/5}(G)} \, dz &\le C \|v_z\|_{L^6(\Gamma_b)} \int_{-h}^0 \|v\|_{L^6(G)}^5 \, dz \\
		&\le \frac19 \big\| \nabla |v_z|^{3/2} \big\|_{L^2(\Omega)}^2 + C \Big( \int_{-h}^0 \|v\|_{L^6(G)}^5 \, dz \Big)^{3/2} \\
		&\le \frac19 \big\| \nabla |v_z|^2 \big\|_{L^2(\Omega)}^2 + C \|1\|_{L^6_z}^{3/2} \|v\|_{L^6(\Omega)}^{15/2} \\
		&\le \frac19 \big\| \nabla |v_z|^2 \big\|_{L^2(\Omega)}^2 + C \|v\|_{H^1(\Omega)}^{3/2} \|v\|_{L^6(\Omega)}^6.
\end{align*}
Since $\alpha^{3/2} \le (\alpha^2 + \alpha)/2$, the contribution of the last term above is the same as that of the first term in the right-hand side of \eref{eq4: proof of main thm for Dirichlet BC case}.
Consequently,
\begin{equation} \label{end of L6 estimate for Dirichlet BC case}
	\frac16 \frac{d}{dt} \|v\|_{L^6(\Omega)}^6 + \big\| |v|^2 |\nabla v| \big\|_{L^2(\Omega)}^2 \le C \|v\|_{H^1(\Omega)}^2 (\|v\|_{L^6(\Omega)}^6 + 1)
		+ \frac19 \big\| \nabla |v_z|^2 \big\|_{L^2(\Omega)}^2.
\end{equation}

Next, in place of \eref{eq2: proof of main prop} for $r = 3$ we have
\begin{equation} \label{eq5: proof of main thm for Dirichlet BC case}
\begin{aligned}
	\frac13 \frac{d}{dt} \|v_z\|_{L^3(\Omega)}^3 + \frac49 \big\| \nabla |v_z|^{3/2} \big\|_{L^2(\Omega)}^2 + \big\| |v_z|^{1/2} |\nabla v_z| \big\|_{L^2(\Omega)}^2
		&\le C \big( |v| |v_z|^{3/2}, |v_z|^{1/2} |\nabla_H v_z| \big)_\Omega + \frac1h \big(\nabla_H \pi, |v_z| v_z \big)_{\Gamma_b} \\
		&\le \frac1{10} \big\| |v_z|^{1/2} |\nabla_H v_z| \big\|_{L^2(\Omega)}^2 + C \big\| |v| |v_z|^{3/2} \big\|_{L^2(\Omega)}^2 \\
		&\qquad	+ \frac1h \big(\nabla_H \pi_1, |v_z| v_z \big)_{\Gamma_b} + \frac1h \big(\nabla_H \pi_2, |v_z| v_z \big)_{\Gamma_b}.
\end{aligned}
\end{equation}
The second term in the right-hand side is estimated, in virtue of $H^{3/8}(\Omega) \hookrightarrow L^{8/3}(\Omega)$ and $H^{1/2}(\Omega) \hookrightarrow L^3(\Omega)$, by
\begin{equation} \label{eq: key of L3 estimate for vz}
\begin{aligned}
	C \big( |v|^2 |v_z|, |v_z|^2 \big)_\Omega &\le \frac1{10} \big\| |v|^2 |v_z| \big\|_{L^2(\Omega)}^2 + C \|v_z\|_{L^4(\Omega)}^4
		= \frac1{10} \big\| |v|^2 |v_z| \big\|_{L^2(\Omega)}^2 + C \big\| |v_z|^{3/2} \big\|_{L^{8/3}(\Omega)}^{8/3} \\
	&\le \frac1{10} \big\| |v|^2 |v_z| \big\|_{L^2(\Omega)}^2 + C \underbrace{ \Big[ \big\| |v_z|^{\frac32} \big\|_{L^2(\Omega)}^{\frac58} \big\| \nabla |v_z|^{\frac32} \big\|_{L^2(\Omega)}^{\frac38} \Big]^{\frac83} }_{ = \|v_z\|_{L^3(\Omega)}^{5/2} \| \nabla |v_z|^{3/2} \|_{L^2(\Omega)} } \\
	&\le \frac1{10} \big\| |v|^2 |v_z| \big\|_{L^2(\Omega)}^2 + \frac19 \big\| \nabla |v_z|^{3/2} \big\|_{L^2(\Omega)}^2 + \underbrace{ C \|v_z\|_{L^3(\Omega)}^4 }_{ \le C \|v_z\|_{L^2(\Omega)}^2 \|\nabla v_z\|_{L^2(\Omega)}^2 } \hspace{-6mm}\cdot\hspace{2mm} \|v_z\|_{L^3(\Omega)}.
\end{aligned}
\end{equation}
Let us address the remaining two pressure terms in \eref{eq5: proof of main thm for Dirichlet BC case}.
The term for $\pi_1$ is bounded, in virtue of the trace inequality $\|f\|_{L^{8/3}(\Gamma_b)} \le C \|f\|_{L^2(\Omega)}^{1/4} \|\nabla f\|_{L^2(\Omega)}^{3/4}$ if $f = 0$ on $\Gamma_u$ (cf.\ \cite[p.\ 1109]{HiKa2016}), by
\begin{align*}
	\|\nabla_H \pi_1\|_{L^2(G)} \big\| |v_z|^2 \big\|_{L^2(\Gamma_b)} &\le C \int_{-h}^0 \big\| |v| |\nabla_H v| \big\|_{L^2(G)} \, dz \, \big\| |v_z|^{\frac32} \big\|_{L^{\frac83}(\Gamma_b)}^{\frac43} \\
		&\le C \|1\|_{L^2_z} \big\| |v| |\nabla_H v| \big\|_{L^2(\Omega)} \, \underbrace{ \Big[ \big\| |v_z|^{\frac32} \big\|_{L^2(\Omega)}^{\frac14} \big\| \nabla |v_z|^{\frac32} \big\|_{L^2(\Omega)}^{\frac34} \Big]^{\frac43} }_{ = \|v_z\|_{L^3(\Omega)}^{1/2} \| \nabla |v_z|^{3/2} \|_{L^2(\Omega)} } \\
		&\le C \big\| |v| |\nabla_H v| \big\|_{L^2(\Omega)}^2 \|v_z\|_{L^3(\Omega)} + \frac19 \big\| \nabla |v_z|^{3/2} \big\|_{L^2(\Omega)}^2.
\end{align*}
The term for $\pi_2$ is bounded by
\begin{align*}
	\|\nabla_H \pi_2\|_{L^3(G)} \|v_z\|_{L^3(\Gamma_b)}^2 &\le C \|v_z\|_{L^3(\Gamma_b)}^3 = C \big\| |v_z|^{3/2} \big\|_{L^2(\Gamma_b)}^2 \le C \big\| |v_z|^{3/2} \big\|_{L^2(\Omega)} \big\| \nabla |v_z|^{3/2} \big\|_{L^2(\Omega)} \\
		&\le \underbrace{ C \|v_z\|_{L^3(\Omega)}^2 }_{ \le C \|v_z\|_{L^2(\Omega)} \|\nabla v_z\|_{L^2(\Omega)} } \hspace{-6mm}\cdot\hspace{2mm} \|v_z\|_{L^3(\Omega)}
			+ \frac19 \big\| \nabla |v_z|^{3/2} \big\|_{L^2(\Omega)}^2.
\end{align*}
Therefore, we deduce that
\begin{equation} \label{end of L3 estimate for Dirichlet BC case}
\begin{aligned}
	&\frac13 \frac{d}{dt} \|v_z\|_{L^3(\Omega)}^3 + \frac19 \big\| \nabla |v_z|^{3/2} \big\|_{L^2(\Omega)}^2 + \frac9{10} \big\| |v_z|^{1/2} |\nabla v_z| \big\|_{L^2(\Omega)}^2 \\
	\le \;& C \big( \big\| |v| |\nabla_H v| \big\|_{L^2(\Omega)}^2 + \|v_z\|_{L^2(\Omega)}^2 \|\nabla v_z\|_{L^2(\Omega)}^2 + \|v_z\|_{L^2(\Omega)}^2 + \|\nabla v_z\|_{L^2(\Omega)}^2 \big) \|v_z\|_{L^3(\Omega)}
		+ \frac1{10} \big\| |v|^2 |v_z| \big\|_{L^2(\Omega)}^2.
\end{aligned}
\end{equation}

We are now in a position to add \eref{end of L6 estimate for Dirichlet BC case} and \eref{end of L3 estimate for Dirichlet BC case} to arrive at
\begin{equation*}
	\frac{dY}{dt} \le C \|v\|_{H^1(\Omega)}^2 Y + C \big( \big\| |v| |\nabla_H v| \big\|_{L^2(\Omega)}^2 + \|v_z\|_{L^2(\Omega)}^2 \|\nabla v_z\|_{L^2(\Omega)}^2 + \|\nabla v_z\|_{L^2(\Omega)}^2 + \|v\|_{H^1(\Omega)}^2 \big) Y^{1/3},
\end{equation*}
where we set $Y(t) := \frac16 \|v(t)\|_{L^6(\Omega)}^6 + \frac13 \|v_z(t)\|_{L^3(\Omega)}^3 + 1$.
After division by $Y^{1/3}$, it follows from the Gronwall inequality that, for all $t \ge 0$,
\begin{equation*}
	Y(t)^{2/3} \le e^{C \|a\|_{L^2(\Omega)}^2} \big( Y(0)^{2/3} + C B_1(a) + C B_1(a)^2 + C\|a\|_{L^2(\Omega)}^2 \big).
\end{equation*}
We thus obtain
\begin{equation*}
	\|v(t)\|_{L^6(\Omega)}^4 + \|v_z(t)\|_{L^3(\Omega)}^2 \le C e^{C \|a\|_{L^2(\Omega)}^2} \big( \|a\|_{L^6(\Omega)}^4 + \|a_z\|_{L^3(\Omega)}^2 + B_1(a) +  B_1(a)^2 \big) =: B_2(a),
\end{equation*}
which involves only single exponential functions.

\subsection{Step 3: estimate for $\nabla_H v \in L^\infty_t L^2(\Omega)$} \label{subsec3.3}
This part is almost the same as in Subsection \ref{subsec2.3}, except that we replace \eref{eq: estimate of wvz with L4}, using $H^{2/3}(G) \hookrightarrow L^6(G)$, by
\begin{equation*}
	\|w v_z\|_{L^2(\Omega)} \le
		C \|w\|_{L^\infty_z L^6_{xy}} \|v_z\|_{L^3(\Omega)} \le C \|\operatorname{div}_H v\|_{L^2_z L^6_{xy}} \|v_z\|_{L^3(\Omega)}
		\le C \|\nabla_H v\|_{L^2(\Omega)}^{1/3} \|\nabla_H v\|_{H^1(\Omega)}^{2/3} \|v_z\|_{L^3(\Omega)},
\end{equation*}
whence
\begin{align*}
	(w v_z, \Delta_H v)_\Omega &\le C \|\nabla_H v\|_{L^2(\Omega)}^2 (\|v_z\|_{L^3(\Omega)}^6 + \|v_z\|_{L^3(\Omega)}^2) + \frac1{10} \|\nabla \nabla_H v\|_{L^2(\Omega)}^2.
\end{align*}
Now, integration of \eref{eq3: proof of main prop} in $t$ as before yields, for all $t \ge 0$,
\begin{align*}
	\frac12 \|\nabla_H v(t)\|_{L^2(\Omega)}^2 + \frac8{10} \int_0^t \|\nabla \nabla_H v\|_{L^2(\Omega)}^2 \, ds
		&\le \frac12 \|\nabla_H a\|_{L^2(\Omega)}^2 + C B_1(a) + C \|a\|_{L^2(\Omega)}^2 (B_2(a)^3 + B_2(a)),
\end{align*}
which together with \eref{eq: L4 estimate for v in Drichlet BC case} completes the desired estimates for $v \in L^\infty(0, T; H^1(\Omega)) \cap L^2(0, T; H^2(\Omega))$.

\section{Proof of \tref{sec1: proof of main thm} for Neumann BC case} \label{sec2: proof of main thm}
Let us return to the Neumann BC case \eref{eq: Neumann BC} and eliminate the factor $e^{\sqrt t}$ from the upper bound of \pref{main prop}.
The cause of this factor is the lower-order term $\|\nabla_H v\|_{L^2(\Omega)} \|v\|_{L^q(\Omega)}^q$ in \eref{eq: lower-order term in Lq estimate}.
To address it we desire a Poincar\'e-type inequality and an embedding $H^1(\Omega) \hookrightarrow L^q(\Omega)$.
For this reason we employ an $L^6(\Omega)$-estimate as in \sref{sec1: proof of main thm} for $\tilde v = v - \bar v$, instead of $v$ itself, which satisfies
\begin{equation} \label{eq: Poincare for tilde v}
	\|\tilde v\|_{L^q(\Omega)} \le C \|v_z\|_{L^q(\Omega)} \quad (q \in [1, \infty]).
\end{equation}
The above inequality follows from integration in $G$ of the one-dimensional Wirtinger inequality
\begin{equation*}
	\int_{-h}^0 |\tilde v(x, y, z)|^q \, dz \le h^q \int_{-h}^0 |v_z(x, y, z)|^q \, dz \quad \text{a.e.} \; (x, y) \in G.
\end{equation*}
We also observe that $\tilde v$ satisfies (cf.\ \cite[(6.4)]{HiKa2016})
\begin{equation*}
	\partial_t \tilde v - \Delta \tilde v + u \cdot \nabla \tilde v = - \tilde v \cdot \nabla_H \bar v + \frac1h \int_{-h}^0 \operatorname{div}_H (\tilde v \otimes \tilde v) \, dz \quad \text{in} \quad \Omega \times (0, T).
\end{equation*}

\subsection{Step 1: estimate for $\tilde v \in L^\infty_t L^q(\Omega) \, (q \in (2, 6])$}
We multiply the above equation by $|\tilde v|^{q-2} \tilde v$ and integrate over $\Omega$ to obtain
\begin{equation*}
	\frac1q \frac{d}{dt} \|\tilde v\|_{L^q(\Omega)}^q + \frac{4(q-2)}{q^2} \big\| \nabla |\tilde v|^{\frac{q}2} \big\|_{L^2(\Omega)}^2 + \big\| |\tilde v|^{\frac{q}2 - 1} |\nabla \tilde v| \big\|_{L^2(\Omega)}^2
		\le C \Big( \int_{-h}^0 |\tilde v| |\nabla_H \tilde v| \, dz, |\tilde v|^{q-1} \Big)_\Omega - (\tilde v \cdot \nabla_H \bar v, |\tilde v|^{q-2} \tilde v)_\Omega.
\end{equation*}
The first term in the right-hand side is bounded, thanks to calculation similar to Subsection \ref{subsec2.1}, by
\begin{equation} \label{eq: treatment of lower-order term}
\begin{aligned}
	&C (\|\nabla_H \tilde v\|_{L^2(\Omega)}^2 + \|\nabla_H \tilde v\|_{L^2(\Omega)}) \|\tilde v\|_{L^q(\Omega)}^q + \frac{q - 2}{q^2} \big\| \nabla_H |\tilde v|^{\frac{q}2} \big\|_{L^2(\Omega)}^2 \\
	\le \; &C \|\nabla_H \tilde v\|_{L^2(\Omega)}^2 (\|\tilde v\|_{L^q(\Omega)}^q + 1) + \frac{q - 2}{q^2} \big\| \nabla_H |\tilde v|^{\frac{q}2} \big\|_{L^2(\Omega)}^2,
\end{aligned}
\end{equation}
where we have used $\|\tilde v\|_{L^q(\Omega)} \le C\|\tilde v\|_{H^1(\Omega)} \le C \|\nabla \tilde v\|_{L^2(\Omega)}$ (recall \eref{eq: Poincare for tilde v}) for $q \le 6$ to address the lower-order term.
The second term is bounded, again by calculation similar to Subsection \ref{subsec2.1}, by
\begin{align*}
	\|\nabla_H \bar v\|_{L^2(G)} \int_{-h}^0 \big\| |\tilde v|^q \big\|_{L^2(G)} \, dz &\le
		C (\|\nabla_H \bar v\|_{L^2(G)}^2 + \|\nabla_H \bar v\|_{L^2(G)}) \|\tilde v\|_{L^q(\Omega)}^q + \frac{q - 2}{q^2} \big\| \nabla_H |\tilde v|^{\frac{q}2} \big\|_{L^2(\Omega)}^2 \\
	&\le C \|\nabla_H \bar v\|_{L^2(G)}^2 (\|\tilde v\|_{L^q(\Omega)}^q + 1) + \frac{q - 2}{q^2} \big\| \nabla_H |\tilde v|^{\frac{q}2} \big\|_{L^2(\Omega)}^2.
\end{align*}
Therefore, it follows that
\begin{equation*}
	\frac1q \frac{d}{dt} \|\tilde v\|_{L^q(\Omega)}^q + \frac{2(q-2)}{q^2} \big\| \nabla |\tilde v|^{\frac{q}2} \big\|_{L^2(\Omega)}^2 + \big\| |\tilde v|^{\frac{q}2 - 1} |\nabla \tilde v| \big\|_{L^2(\Omega)}^2
		\le C \|\nabla_H v\|_{L^2(\Omega)}^2 (\|\tilde v\|_{L^q(\Omega)}^q + 1),
\end{equation*}
which yields, by the Gronwall inequality, for all $t \ge 0$,
\begin{equation*}
	\|\tilde v(t)\|_{L^q(\Omega)}^q + \int_0^t \big\| |\tilde v|^{\frac{q}2 - 1} |\nabla \tilde v| \big\|_{L^2(\Omega)}^2 \, ds
		\le C \big( \|\tilde a\|_{L^q(\Omega)}^q + \|a\|_{L^2(\Omega)}^2 \big)  \big( 1+ \|a\|_{L^2(\Omega)}^2 e^{C \|a\|_{L^2(\Omega)}^2} \big) =: E_1(q, a).
\end{equation*}

\subsection{Step 2: estimates for $v_z \in L^\infty_t L^r(\Omega) \, (r = 2, 3)$}
We find from \eref{eq2: proof of main prop} that, since $\operatorname{div}_H v = \operatorname{div}_H \tilde v$,
\begin{equation} \label{eq1: Step 2 of main thm for Neumann}
\begin{aligned}
	&\frac1r \frac{d}{dt} \|v_z\|_{L^r(\Omega)}^r + \frac{4(r-2)}{r^2} \big\| \nabla |v_z|^{\frac{r}2} \big\|_{L^2(\Omega)}^2 + \big\| |v_z|^{\frac{r}2 - 1} |\nabla v_z| \big\|_{L^2(\Omega)}^2 \\
		= \; &-(v_z \cdot \nabla_H v, |v_z|^{r-2} v_z)_\Omega + (\operatorname{div}_H v, |v_z|^r)_\Omega \\
		= \; &-(v_z \cdot \nabla_H \bar v, |v_z|^{r-2} v_z)_\Omega - (v_z \cdot \nabla_H \tilde v, |v_z|^{r-2} v_z)_\Omega + (\operatorname{div}_H \tilde v, |v_z|^r)_\Omega.
\end{aligned}
\end{equation}
If $r = 2$, then the first term and the last two terms in the right-hand side are bounded by
\begin{align*}
	\|\nabla_H \bar v\|_{L^2(G)} \int_{-h}^0 \|v_z\|_{L^4(G)}^2 \, dz &\le \|\nabla_H \bar v\|_{L^2(G)} \int_{-h}^0 \|v_z\|_{L^2(G)} \|v_z\|_{H^1(G)} \, dz
		\le \|\nabla_H \bar v\|_{L^2(G)} \|v_z\|_{L^2(\Omega)} \|\nabla v_z\|_{L^2(\Omega)} \\
	&\le C \|\nabla_H \bar v\|_{L^2(G)}^2 \|v_z\|_{L^2(\Omega)}^2 + \frac1{10} \|\nabla v_z\|_{L^2(\Omega)}^2,
\end{align*}
and by
\begin{equation*}
	C \big\| |\tilde v| |v_z| \big\|_{L^2(\Omega)}^2 + \frac1{10} \|\nabla_H v_z\|_{L^2(\Omega)}^2,
\end{equation*}
respectively. Therefore, we have
\begin{equation*}
	\frac12 \frac{d}{dt} \|v_z\|_{L^2(\Omega)}^2 + \frac8{10} \|\nabla v_z\|_{L^2(\Omega)}^2 \le C \|\nabla v\|_{L^2(\Omega)}^2 \|v_z\|_{L^2(\Omega)}^2 + C \big\| |\tilde v| \, |v_z| \big\|_{L^2(\Omega)}^2,
\end{equation*}
which implies (note that $v_z = \tilde v_z$), for all $t \ge 0$,
\begin{equation*}
	\|v_z(t)\|_{L^2(\Omega)}^2 + \int_0^t \|\nabla v_z\|_{L^2(\Omega)}^2 \, ds
		\le C \big( \|a_z\|_{L^2(\Omega)}^2 + E_1(4, a) \big) \big( 1 + \|a\|_{L^2(\Omega)}^2 e^{C \|a\|_{L^2(\Omega)}^2} \big) =: E_2(a),
\end{equation*}
by the Gronwall inequality.

If $r = 3$, then the first term and the last two terms of \eref{eq1: Step 2 of main thm for Neumann} are bounded by
\begin{align*}
	\|\nabla_H \bar v\|_{L^2(G)} \int_{-h}^0 \big\| |v_z|^{\frac32} \big\|_{L^4(G)}^2 \, dz
		&\le C\|\nabla v\|_{L^2(\Omega)} \int_{-h}^0 \big\| |v_z|^{\frac32} \big\|_{L^2(G)} \big\| |v_z|^{\frac32} \big\|_{H^1(G)} \, dz \\
	&\le C\|\nabla v\|_{L^2(\Omega)} \|v_z\|_{L^3(\Omega)}^{\frac32} \big\| \nabla |v_z|^{\frac32} \big\|_{L^2(\Omega)} \\
	&\le C\|\nabla v\|_{L^2(\Omega)}^2 \|v_z\|_{L^3(\Omega)}^3 + \frac19 \big\| \nabla |v_z|^{\frac32} \big\|_{L^2(\Omega)}^2,
\end{align*}
and by (recall \eref{eq: key of L3 estimate for vz} for similar calculation)
\begin{align*}
	C \big( |\tilde v| |v_z|^{3/2}, |v_z|^{1/2} |\nabla_H v_z| \big)_\Omega &\le C \big\| |\tilde v| |v_z|^{3/2} \|_{L^2(\Omega)}^2 + \frac1{10} \big\| |v_z|^{1/2} |\nabla_H v_z| \big\|_{L^2(\Omega)}^2 \\
		&\le C \|v_z\|_{L^4(\Omega)}^4 + \big\| |\tilde v|^2 |v_z| \big\|_{L^2(\Omega)}^2 + \frac1{10} \big\| |v_z|^{1/2} |\nabla_H v_z| \big\|_{L^2(\Omega)}^2 \\
		&\le C \|v_z\|_{L^2(\Omega)}^2 \|\nabla v_z\|_{L^2(\Omega)}^2 \|v_z\|_{L^3(\Omega)} + \frac19 \big\| \nabla |v_z|^{3/2} \big\|_{L^2(\Omega)}^2 \\
		&\qquad + \big\| |\tilde v|^2 |v_z| \big\|_{L^2(\Omega)}^2 + \frac1{10} \big\| |v_z|^{1/2} |\nabla_H v_z| \big\|_{L^2(\Omega)}^2,
\end{align*}
respectively. Therefore, we have
\begin{align*}
	&\frac13 \frac{d}{dt} \|v_z\|_{L^3(\Omega)}^3 + \frac29 \big\| \nabla |v_z|^{3/2} \big\|_{L^2(\Omega)}^2 + \frac8{10} \big\| |v_z|^{1/2} |\nabla v_z| \big\|_{L^2(\Omega)}^2 \\
		\le \; &C\|\nabla v\|_{L^2(\Omega)}^2 \|v_z\|_{L^3(\Omega)}^3 + C \|v_z\|_{L^2(\Omega)}^2 \|\nabla v_z\|_{L^2(\Omega)}^2 \|v_z\|_{L^3(\Omega)} + \big\| |\tilde v|^2 |v_z| \big\|_{L^2(\Omega)}^2,
\end{align*}
in particular,
\begin{equation*}
	\frac{d}{dt} (\|v_z\|_{L^3(\Omega)}^3 + 1) \le C\|\nabla v\|_{L^2(\Omega)}^2 (\|v_z\|_{L^3(\Omega)}^3 + 1) + 
		C \Big( \|v_z\|_{L^2(\Omega)}^2 \|\nabla v_z\|_{L^2(\Omega)}^2 + \big\| |\tilde v|^2 |v_z| \big\|_{L^2(\Omega)}^2 \Big) \big( \|v_z\|_{L^3(\Omega)}^3 + 1 \big)^{1/3}.
\end{equation*}
Dividing the both sides by $\big( \|v_z\|_{L^3(\Omega)}^3 + 1 \big)^{1/3}$ and applying the Gronwall inequality, we obtain, for all $t \ge 0$,
\begin{equation*}
	\|v_z(t)\|_{L^3(\Omega)}^2 \le C e^{C \|a\|_{L^2(\Omega)}^2} \big( \|a_z\|_{L^3(\Omega)}^2 + E_2(a)^2 + E_1(6, a) \big) =: E_3(a).
\end{equation*}

\subsection{Step 3: estimate for $\nabla_H v \in L^\infty_t L^2(\Omega)$}
In place of \eref{eq3: proof of main prop} we have
\begin{align*}
	&\frac12 \frac{d}{dt} \|\nabla_H v\|_{L^2(\Omega)}^2 + \|\nabla \nabla_H v\|_{L^2(\Omega)}^2
		= (\bar v \cdot \nabla_H v + \tilde v \cdot \nabla_H \bar v, \Delta_H v)_\Omega + (\tilde v \cdot \nabla_H \tilde v + w v_z, \Delta_H v)_\Omega \\
	\le \; & (\bar v \cdot \nabla_H v + \tilde v \cdot \nabla_H \bar v, \Delta_H v)_\Omega
		 + C \big\| |\tilde v| |\nabla_H \tilde v| \big\|_{L^2(\Omega)}^2 + C \|\nabla_H v\|_{L^2(\Omega)}^2 (\|v_z\|_{L^3(\Omega)}^6 + \|v_z\|_{L^3(\Omega)}^2) + \frac2{10} \|\nabla \nabla_H v\|_{L^2(\Omega)}^2,
\end{align*}
where we have addressed $\|w v_z\|_{L^2(\Omega)}$ in the same way as in Subsection \ref{subsec3.3}.
From $\operatorname{div}_H \bar v = 0$ and $\alpha^3 \le (\alpha^4 + \alpha^2)/2$, we find that
\begin{align*}
	(\bar v \cdot \nabla_H v, \Delta_H v)_\Omega &= -\int_\Omega (\nabla_H \bar v) (\nabla_H v) (\nabla_H v)^\top - \underbrace{ \int_\Omega \frac12 (\bar v \cdot \nabla_H) |\nabla_H v|^2 }_{= 0} \\
		&\le C \|\nabla_H \bar v\|_{L^2(G)} \int_{-h}^0 \|\nabla_H v\|_{L^4(G)}^2 \, dz \le C \|\nabla_H \bar v\|_{L^2(G)} \|\nabla_H v\|_{L^2(\Omega)} \|\nabla_H v\|_{L^2_zH^1_{xy}} \\
		&\le C \|\nabla_H v\|_{L^2(\Omega)}^4 + \frac1{10} \|\nabla_H^2 v\|_{L^2(\Omega)}^2 + C \|\nabla_H v\|_{L^2(\Omega)}^2.
\end{align*}
Finally, observe that
\begin{align*}
	(\tilde v \cdot \nabla_H \bar v, \Delta_H v)_\Omega &\le \|\tilde v\|_{L^4(\Omega)} \|\nabla_H \bar v\|_{L^4(G)} \|\Delta_H v\|_{L^2(\Omega)} \\
		&\le C \|\tilde v\|_{L^4(\Omega)} \|\nabla_H \bar v\|_{L^2(G)}^{1/2} \|\nabla_H \bar v\|_{H^1(G)}^{1/2} \|\Delta_H v\|_{L^2(\Omega)} \\
		&\le C \|\nabla_H \bar v\|_{L^2(G)}^2 (\|\tilde v\|_{L^4(\Omega)}^4 + \|\tilde v\|_{L^4(\Omega)}^2) + \frac2{10} \|\nabla_H^2 v\|_{L^2(\Omega)}^2.
\end{align*}

Collecting the estimates above, together with $\alpha^2 \le (\alpha^6 + 2)/3$ and $\alpha^2 \le (\alpha^4 + 1)/2$, leads to
\begin{align*}
	\frac12 \frac{d}{dt} \|\nabla_H v\|_{L^2(\Omega)}^2 + \frac5{10} \|\nabla \nabla_H v\|_{L^2(\Omega)}^2
		&\le C \|\nabla_H v\|_{L^2(\Omega)}^4 + C \big\| |\tilde v| |\nabla_H \tilde v| \big\|_{L^2(\Omega)}^2 \\
		&\qquad	+ C \|\nabla_H v\|_{L^2(\Omega)}^2 (\|v_z\|_{L^3(\Omega)}^6 + \|\tilde v\|_{L^4(\Omega)}^4 + 1).
\end{align*}
This combined with the Gronwall inequality implies that, for all $t \ge 0$, 
\begin{align*}
	&\|\nabla_H v(t)\|_{L^2(\Omega)}^2 + \int_0^t \|\nabla \nabla_H v\|_{L^2(\Omega)}^2 \, ds \\
	\le \; & C (1 + \|a\|_{L^2(\Omega)}^2 e^{C \|a\|_{L^2(\Omega)}^2})
		\big( \|\nabla_H a\|_{L^2(\Omega)}^2 + E_1(4, a) + \|a\|_{L^2(\Omega)}^2 (E_3(a)^3 + E_1(4, a) + 1) \big),
\end{align*}
which concludes the desired estimate for $v \in L^\infty(0, T; H^1(\Omega)) \cap L^2(0, T; H^2(\Omega))$.


\end{document}